\documentclass{article}
\usepackage{amssymb}
\usepackage{amsmath}

\let\<\langle
\let\>\rangle
\font\bigcmsy=cmsy10.pk scaled 2000
\def\bigtimes{\mathop{\,\vrule width0pt depth2pt height8pt
            \smash{\lower2pt\hbox{\bigcmsy\char'002}}\,}\limits}

\begin{document}

\begin{center}
\Large{\textbf{Order separability.}}
\end{center}

\begin{center}
\textbf{Vladimir V. Yedynak}
\end{center}

\begin{abstract}
This paper is devoted to the investigation of the property of order separability for free products of groups.


\end{abstract}

\section{Introduction.}

Definition. A group $G$ is called order separable if for each elements $u$ and $v$ of $G$ such that $u$ is conjugate to neither $v$ nor $v$ inverse there exists a homomorphism $\varphi$ of $G$ onto a finite group such that the orders of $\varphi(u)$ and $\varphi(v)$ are different.

In [1] it was proved that free groups are order separable. In this work we prove that this property is inherited by free products:

\textbf{Theorem 1.} \textsl{The group $G=A\ast B$ is order separable if and only if $A$ and $B$ are order separable.}

Note that the property of order separability for free groups was generalized in [2] where it was proved that free groups are actually omnipotent.

\section{Notations and definitions.}

Investigate the graph $\Gamma$ satisfying the following properties:

1) $\Gamma$ is an oriented graph whose positively oriented edges are labelled by elements of groups $A$ and $B$ so that for each vertex $p$ of $\Gamma$ and for each $a\in A$ and $b\in B$ there exist exactly one edge with label $a$ and exactly one edge with label $b$ ending at $p$ and there exist exactly one edge with label $a$ and exactly one edge with label $b$ starting at $p$;

2) for each vertex $p$ of $\Gamma$ we define the subgraph $A(p)$ of the graph $\Gamma$ as the maximal connected graph which contains $p$ and whose positively oriented edges are labelled by the elements of $A$; it is required that $A(p)$ is the Cayley graph of $A$ with the set of generators $\{A\}$. The graph $B(p)$ is defined analogically.

We shall use the following notations. The symbols Lab $(e), \alpha(f), \omega(f), \alpha(S), \omega(S)$ will denote correspondingly the label of the positively oriented edge $e$, the beginning and the end of the edge $f$ and the beginning and the end of the path $S$. Having a path $S=e_1...e_k$ we define its label Lab $(S)=$ Lab $(e_1)...$ Lab $(e_n)$.

Definition 1. Consider the graph $\Gamma$ satisfying the properties 1), 2) and the cyclically reduced element $u\in A\ast B$. The closed path $S=e_1...e_n$ is called $u$-cycle if Lab $(e_{il+1}...e_{il+l})=u$ where $l$ is the length of the element $u$, $k$ is an arbitrary natural number and subscripts are modulo $n$.

If a label of the $u$-cycle $S$ is $u^k$ then we shall say that the length of the $u$-cycle $S$ equals $k$.

The group $A\ast B$ acts on the right on the set of vertices of the graph $\Gamma$ by the following way. Consider the vertex $p$ of $\Gamma$ and the elements $c\in(A\cup B)\setminus\{1\}$. Then according to the property 1) there exist the edge $u$ with label $c$ starting at $p$ and the edge $v$ ending at $p$ and the labels of $u$ and $v$ coincide with $c$. Then we put $p\circ c=\omega(u), p\circ c^{-1}=\alpha(v)$.

Definition 2. We say that the cycle $S=e_1...e_n$ of a graph $\Gamma$ with properties 1), 2) does not have near edges if there are no distinct edges of $S$ belonging to one subgraph $A(p)$ or $B(p)$ for some $p$.

Definition 3. A group $G$ is called subgroup separable if each finitely generated subgroup of $G$ coincides with the intersection of finite index subgroups of $G$.

In [3] the following theorem was proved.

\textbf{Theorem 2.} \textsl{The class of subgroup separable groups is closed with respect to the operation of the free product of groups.}

\textbf{Corollary}. \textsl{The free product of finite groups is cyclic subgroup separable.}

\section{Proof of theorem 1.}

If $A\ast B$ is order separable then it is obvious that $A$ and $B$ are order separable. Consider order separable groups $A$ and $B$ and prove that $A\ast B$ is order separable. Put $G=A\ast B$. Consider cyclically reduced elements $u$ and $v$ of $G$ such that $u$ is not conjugate to $v^{\pm1}$. If $u$ and $v$ belong to free factors then we use the natural homomorphism of $G$ onto $A$ or $B$ and use the order separability of free factors.

Suppose that $u\notin A\cup B$. Consider the case when $u$ and $v$ belong to the Cartesian subgroup $C=\<[a,b]| a\in A, b\in B\>$ and do not equal to unit. Consider that for each homomorphism of $G$ onto a finite group the images of $u$ and $v$ have equal orders. It is possible to consider that the normal forms for $u$ and $v$ have the following presentations: $u=a_1b_1...a_nb_n, v=a_1'b_1'...a_m'b_m'$. Since order separability involves residual finiteness we may deduce that there exists a homomorphism of $G$ onto a group $A_1\ast B_1$ such that $A_1$ is the image of $A$ and $B_1$ is the image of $B$ besides $a_i, b_i, a_j', b_j'$ have nonunit images and each element presented as $a_ia_k, b_ib_k, a_j'a_l', b_j'b_l', a_ia_l', b_ib_l'$ which differs from unit has a nonunit image too, $i, k=1,..., n, j, l=1,..., m$. Thereby we may consider that the groups $A$ and $B$ are finite. For each number $n=0, 1, 2,...$ construct the graph  $\Gamma_n$ with properties 1), 2) which satisfies also the following properties:

3) the length of each $u$-cycle divides the length of a maximal $u$-cycle; the same is true for $v$-cycles;

4) in $\Gamma_n$ (when $n>0$) there exists the path $R_n$ of length $n$ which is contained in a maximal $u$-cycle and in all maximal $v$-cycles;

5) all $u$- and $v$-cycle of $\Gamma_n$ have no near edges;

6) the length of a maximal $u$-cycle coincides with the length of a maximal $v$-cycle

The construction of $\Gamma_0$. Due to the corollary there exists the homomorphism $\varphi$ of $G$ onto a finite group such that the elements  $ax$ and $bz$ which are not conjugate to elements from $\<u\>$ and $\<v\>$ have nonunit images where $x$ and $z$ are the subwords of words $u^k, v^k$, $k=0,1,2,..., a\in A, b\in B$. We may also consider that $u$ and $v$ do not belong to the kernel of $\varphi$. Then we may take the Cayley graph of $\varphi(G)$ with the generating set $\varphi(A\cup B)$ in the capacity of $\Gamma_0$ (labels $\varphi(a), \varphi(b)$ are identified with $a$ and $b$ correspondingly). Conditions 1), 2) and 3) are held because of the definition of the Cayley graph; conditions 4), 5) are held due to the properties of the homomorphism; the property 6) is true by the supposition about the orders of images of $u$ and $v$.

The construction of $\Gamma_{n+1}$ from $\Gamma_n$. Let $t$ be the length of the maximal $u$-cycle in $\Gamma_n$. Consider $t$ copies of $\Gamma_n: \Delta_1,..., \Delta_t$. Put $q_k=\omega(R_{n, k})$ where $R_{n, k}$ is the path in $\Delta_k$ corresponding to the path $R_n$ of $\Gamma_n$, $p_k$ is the vertex following after $q_k$ on the maximal $u$-cycle passing through $R_{n, k}$ (it is supposed that $p_k$ does not belong to $R_{n, k}$ and vertices $p_k$ in graphs $\Delta_k$ and chosen maximal $u$-cycles passing through $R_{n, k}$ correspond to each other). If $n=0$ then $q_k$ ia an arbitrary vertex and the edge $(q_k, p_k)$ belongs to a $u$-cycle. Otherwise $(q_k, p_k)$ is the edge connecting $p_k$ and $q_k$ and belonging to the chosen maximal $u$-cycle. Consider that Lab$(q_k, p_k)\in A$. In order to construct the graph $K_{n, 1}$ from $\Delta_1,..., \Delta_t$ we delete all edges from $A(q_k)$ which are incident to $q_k$. Let $s$ be an arbitrary vertex of the subgraph $A(q)$ of the graph $\Gamma_n$ which differs from $q$ and $s$ is connected with $q$ be the edge $e\in A(q)$. The vertex $s_k\in\Delta_k$ corresponds to the vertex $s$. Connect the vertex $q_k$ by the edge with the vertex $s_{k+1}$ (if $k=t$ we consider that $k+1=1$). The label of this new edge $f_k$ equals Lab $(e)$. Also if $\alpha(e)=q$ in $\Gamma_n$ then $\alpha(f_k)=q_k$; if $\omega(e)=q$ in $\Gamma_n$ then $\omega(f_k)=q_k$. Put $S_n=R_{n, 1}\cup d_1$ where $d_1$ is the edge which is appended instead of the edge $(q_1, p_1)$. The graph $K_{n, 1}$ satisfies properties 1), 2). The property 3) is fulfilled since the lengths of each $u$- or $v$-cycle either does not change or becomes $t$ times greater than it was. The condition 5) is true because otherwise it is not held for $\Gamma_n$. If in the graph $K_{n, 1}$ all maximal $v$-cycles pass through $S_n$ then we put $\Gamma_{n+1}=K_{n, 1}, R_{n+1}=S_n$. Otherwise consider $t^2$ copies of the graph $K_{n, 1}: \Omega_1,..., \Omega_{t^2}$. Let $r_k$ be the vertex next to $\omega(S_{n, k})$ on the maximal $u$-cycle passing through $S_{n, k}$ where $S_{n, k}$ is a path of $\Omega_k$ corresponding to $S_n$ in the graph $K_{n, 1}$ ($r_k$ does not belong to $S_{n, k}$). Vertices $r_k$ and maximal $u$-cycles passing through $S_{n, k}$ correspond to each other in $\Omega_k$. Construct the graph $K_{n, 2}$ from $\Omega_1,..., \Omega_{t^2}$ the same way as the graph $K_{n, 1}$ is constructed from $\Delta_1,..., \Delta_t$ but we consider $\omega(S_{n, k}), r_k, B(\omega(S_{n, k}))$ instead of vertices $q_k, p_k$ and the subgraph $A(q_k)$ correspondingly. The graph $K_{n, 2}$ satisfies the properties 1), 2), 3), 5). This is established in similar way as for the graph $K_{n, 1}$. The path $S_{n, 1}$ is contained in a maximal $u$-cycle and in all maximal $v$-cycles of the graph $K_{n, 2}$ because of the property 5) and 6) for $K_{n, 2}$. Thus $\Gamma_{n+1}=K_{n, 2}, R_{n+1}=S_{n, 1}$.

Since $n$ is an arbitrary natural number then conjugating $u$ and $v$ we may consider that $u=w^k, v=w^l$. It is possible to consider that $(k, l)=1$, that is $k$ and $l$ are coprime. Hence $w\in C$. Suppose that $|k|>1$. Then there exists a prime number $p$ such that $p\mid k, p\nmid l$. There exists a homomorphism $\psi$ of $C$ onto a finite $p$-group $P$ such that $w$ has a nonunit image [4]. Put $N=$ ker $\psi$. Then the group $N'=\cap_{g\in G}g^{-1}Ng$ is a finite index normal divisor of $G$. Besides in the quotient-group $G/N'$ the image of $w$ has the order  which equals the nonzero power of $p$. Denote by $\psi_1$ the natural homomorphism of $G$ onto a finite group $G/N'$. Because of the conditions on the order of $\psi_1(w)$ we conclude that the order of elements $\psi_1(u), \psi_1(v)$ are different. So $|k|=|l|=1$ and this involves the violation.

Consider now the case when $u$ and $v$ belong to $A\ast B\setminus(C\cup A\cup B)$. It was shown in [5] that this condition involves that $u$ and $v$ have infinite orders. Since the groups $A$ and $B$ are finite there exists the natural number $q$ such that $u^q$ and $v^q$ belong to $C$. Besides since $u$ and $v$ are cyclically reduced and $u$ and $v$ do not belong to $A\ast V\setminus(C\cup A\cup B)$ and due to the conjugacy theorem for free products [5] we deduce that $u^n$ is not conjugate to $v^{\pm n}$. If the orders of images of $u^q$ and $v^q$ are different after some homomorphism then the orders of images of $u$ and $v$ are also different after the same homomorphism. The case when $u\notin A\ast B\setminus\{\cup_{g\in G}g^{-1}(A\cup B)g\}$ and $v\in h^{-1}(A\cup B)h, h\in A\ast B$ can be solved with the usage of the residual finiteness of $G$.

Theorem 1 is proved.

\begin{center}
\large{Acknowledgements.}
\end{center}

The author thanks A. A. Klyachko for setting the problem and valuable comments.

\begin{center}
\large{References.}
\end{center}

1. \textsl{Klyachko A. A.} Equations over groups, quasivarieties,
and a residual property of a free group // J. Group Theory. 1999.
\textbf{2}. 319--327.

2. \textsl{Wise, Daniel T.} Subgroup separability of graphs of free groups with cyclic edge groups.
Q. J. Math. 51, No.1, 107-129 (2000). [ISSN 0033-5606; ISSN 1464-3847]

3. \textsl{Romanovskii N. S.} On the residual finiteness of free products with respect to membership. // Izv. AN SSSR. Ser.
matem., 1969, \textbf{33}, 1324-1329.

4. Kargapolov, M. I., Merzlyakov, Yu. I. (1977). Foundations of group
theory. \textsl{Nauka}.

5. Lyndon, R. C., Schupp, P. E. (1977). Combinatorial group
theory. \textsl{Springer-Verlag}.

\end{document}